\documentclass[11pt,a4paper]{amsart}
\pdfoutput=1

\usepackage[utf8]{inputenc}

\usepackage{hyperref,graphicx,amssymb}

\newcommand{\CC}{\mathbb{C}}

\begin{document}

\title{Belyi map for the sporadic group $\text{J}_\mathbf{2}$}

\author{Dominik Barth}
\author{Andreas Wenz}

\address{Institute of Mathematics\\ University of Würzburg \\ Emil-Fischer-Straße 30 \\ 97074 Würzburg, Germany}
\email{dominik.barth@mathematik.uni-wuerzburg.de}
\email{andreas.wenz@mathematik.uni-wuerzburg.de}



\begin{abstract}
We present a genus $0$ Belyi map for the sporadic Janko group $\text{J}_2$ of degree $280$. As a consequence we obtain a polynomial having $\text{Aut}(\text{J}_2)$ as a Galois group over $K(t)$ where $K$ is a number field of degree 10.
\end{abstract}
\maketitle

\section{Introduction}

Using the method described in \cite{Barth2017} we compute a Belyi map of degree 280 with monodromy group isomorphic to the sporadic Janko group $\text{J}_2$.
A Belyi map of degree 100 with monodromy group isomorphic to $\text{J}_2$ has already been computed by Monien, see \cite{Monien2017}.

In this paper we work with the permutation triple $(x,y,z:=(xy)^{-1}) \in S_{280}^3$ given in the ancillary data. It is of type
\begin{center}
\renewcommand{\arraystretch}{1.3}
\begin{tabular}{c|c|c|c}
& $x$ & $y$ & $z = (xy)^{-1}$ \\  \hline
cycle structure &  $2^{134}.1^{12}$ & $7^{40}$  & $3^{92}.1^4$ \\
\end{tabular}
\end{center}
and has the following properties:
\begin{itemize}
\item $x,y$ generate the sporadic Janko group $\text{J}_2$.
\item $(x,y,z)$ is of genus $0$.
\item The permutations $x,y,z$ each lie in rational conjugacy classes.
\item $(x,y,z)$ has passport size $10$.
\end{itemize}

\section{Results and Verification}

The computed Belyi map $$f = \frac{p}{q} = 1 + \frac{r}{q},$$ defined over the number field $K = \mathbb{Q}(\alpha)$ where
$$
\alpha^{10} - \alpha^9 - 4\alpha^8 - 2 \alpha^7 + 65\alpha^6 + 27\alpha^5 + 11\alpha^4 - 89\alpha^3 + 25\alpha^2 - 13\alpha + 1
=0
$$	
can be found in the ancillary file and its reduction modulo a prime ideal lying over 283 is given in the appendix. Note that Monien's Belyi map of degree 100 corresponding to $\text{J}_2$ is defined over the same number field although it is described by a different polynomial, see Theorem 1 in \cite{Monien2017}. 

With the Riemann-Hurwitz genus formula and the factorizations of $p,q$ and $r$ one can easily verify that $f:\mathbb{P}^1(\CC) \rightarrow \mathbb{P}^1(\CC)$ is indeed a Belyi map, i.e.\;a three branch point covering, where the ramification over $0$, $1$ and $\infty$ corresponds to the cycle structures of $x$, $y$ and $z$, respectively.
In the following we will show that
\begin{itemize}
\item the geometric monodromy group $G:= \text{Gal}(p(X)-tq(X)\in \mathbb{C}(t)[X])$ is isomorphic to $\text{J}_2$ and
\item the arithmetic monodromy group 
$A:= \text{Gal}(p(X)-tq(X)\in K(t)[X])$
is isomorphic to $\text{Aut}(\text{J}_2)$.
\end{itemize}

Let $\mathcal{O}_K$ be the ring of integers in $K$, and $\mathfrak{p} := (283,167 + \alpha)\mathcal{O}_K$. Then $\mathfrak{p}$ is a prime ideal in $\mathcal{O}_K$ lying over $283$. Since all coefficients of $p$ and $q$ are contained in the localization of $\mathcal{O}_K$ at $\mathfrak{p}$ we can reduce them modulo $\mathfrak{p}$ to obtain polynomials $\bar{p},\bar{q}\in\mathbb{F}_{283}[X]$ leading us to study 
$$
A_{\mathbb{F}}:= \text{Gal}(\bar{p}(X)-t\bar{q}(X)\in \mathbb{F}_{283}(t)[X]).$$
Computing the irreducible factors of $\bar{p}(t)\bar{q}(X)-\bar{p}(X)\bar{q}(t)\in \mathbb{F}_{283}(t)[X]$ enables us to determine the subdegrees of $A_\mathbb{F}$ which turn out to be 1, 36, 108 and 135. Because there is no subset of the subdegrees containing 1 adding up to a nontrivial divisor of the permutation degree 280, $A_\mathbb{F}$ is primitive. According to the Magma database \cite{Magma} for finite primitive groups and the fact that the discriminant of $\bar{p}-t\bar{q}$ is not a square in $\mathbb{F}_{283}(t)$ only one possibility remains: $A_{\mathbb{F}} =  \text{Aut}(\text{J}_2)$. 
Thanks to Dedekind reduction \cite[VII, Theorem 2.9]{Lang2002} this implies: $A$ is primitive and $\text{Aut}(\text{J}_2)$ is a subgroup of $A$, therefore $A \in \lbrace  \text{Aut}(\text{J}_2), S_{280} \rbrace$. Since $p(t)q(X)-p(X)q(t)\in K(t)[X]$ has a divisor of degree $36$, see ancillary file, $A$ is not $2$-transitive, thus $A = \text{Aut}(\text{J}_2)$.

Taking into account that $G$ is normal in $A$ and $\text{J}_2$ is simple, we also get
$G = \text{J}_2$ or $G = \text{Aut}(\text{J}_2)$. Because permutations with the same cycle structures as $x,y$ and $z$ generate an even permutation group, we end up with $G=\text{J}_2$.

\appendix
\section*{Appendix: Computed Data}
Sticking to the notation of the previous section we will present the computed Belyi map $f = \frac{p}{q} = 1+ \frac{r}{q}\in K(X)$ reduced by the prime ideal $\mathfrak{p}= (283,167+\alpha)\mathcal{O}_K$. The resulting polynomials $\bar{r}$ and $\bar{q}$ are                                                                                         
\begin{align*}
\bar{r}(X) = \;& 157 \cdot (X + 212)^7\\                                                                                                       
&\cdot(X^2 + 134X + 135)^7\\                                                                                                                                                                                                    
&\cdot(X^3 + 66X^2 + 49X + 21)^7\\                                                                                                                                                                                             
&\cdot(X^3 + 70X^2 + 18X + 129)^7\\                                                                                                                                                                                            
&\cdot(X^3 + 75X^2 + 71X + 119)^7\\                                                                                                                                                                                            
&\cdot(X^3 + 220X^2 + 76X + 2)^7\\                                                                                                                                                                                             
&\cdot(X^6 + 141X^5 + 44X^4 + 222X^3 + 128X^2 + 135X + 240)^7\\                                                                                                                                                             
&\cdot(X^6 + 153X^5 + 85X^4 + 215X^3 + 247X^2 + 177X + 94)^7\\                                                                                                                                                              
&\cdot(X^6 + 235X^5 + 142X^4 + 83X^3 + 3X^2 + 41X + 186)^7\\                                                                                                       
&\cdot(X^6 + 267X^5 + 199X^4 + 191X^3 + 34X^2 + 123X + 279)^7
\end{align*}
and
\vspace{-2mm}
\begin{align*}
\bar{q}(X) = \;& (X^4 + 185X^3 + 88X^2 + 11X + 93)\\                                                                            
&\cdot (X^2 + 147X + 261)^3  \\                                                                  
&\cdot (X^2 + 240X + 211)^3  \\                                                                                            
&\cdot (X^4 + X^3 + 258X^2 + 211X + 120)^3  \\                                                                             
&\cdot (X^4 + 5X^3 + 138X^2 + 41X + 206)^3  \\                                                                            
&\cdot (X^4 + 12X^3 + 264X^2 + 18X + 202)^3  \\                                                                           
&\cdot (X^4 + 41X^3 + 55X^2 + 165X + 210)^3  \\                                                                           
&\cdot (X^4 + 60X^3 + 198X^2 + 112X + 213)^3  \\                                                                          
&\cdot (X^4 + 88X^3 + 179X^2 + 79X + 27)^3  \\
&\cdot (X^4 + 90X^3 + 232X^2 + 231X + 36)^3  \\                                                                           
&\cdot (X^4 + 104X^3 + 64X^2 + 124X + 138)^3  \\                                                                          
&\cdot (X^4 + 109X^3 + 78X^2 + 111X + 275)^3  \\                                                                          
&\cdot (X^4 + 110X^3 + 275X^2 + 62X + 76)^3  \\
&\cdot (X^4 + 121X^3 + 277X^2 + 62X + 234)^3  \\                                                                          
&\cdot (X^4 + 145X^3 + 66X^2 + 211X + 79)^3  \\                                                                           
&\cdot (X^4 + 150X^3 + 210X^2 + 78X + 72)^3  \\                                                                          
&\cdot (X^4 + 153X^3 + 176X^2 + 67X + 75)^3  \\                                                                           
&\cdot (X^4 + 156X^3 + 14X^2 + 140X + 170)^3  \\                                                                          
&\cdot (X^4 + 160X^3 + 107X^2 + 72X + 223)^3  \\                                                                          
&\cdot (X^4 + 180X^3 + 77X^2 + 38X + 49)^3  \\                                                                            
&\cdot (X^4 + 184X^3 + 156X^2 + 43X + 194)^3  \\                                                                          
&\cdot (X^4 + 189X^3 + 76X^2 + 125X + 145)^3  \\                                                                          
&\cdot (X^4 + 201X^3 + 162X^2 + 88X + 127)^3  \\                                                                          
&\cdot (X^4 + 220X^3 + 276X^2 + 186X + 9)^3  \\                                                                           
&\cdot (X^4 + 234X^3 + 23X^2 + 219X + 118)^3 .
\end{align*}

\end{document}